\documentclass[a4paper,11pt]{amsart}
\usepackage[cp1251]{inputenc}
\usepackage[english]{babel}
\usepackage{amsmath,amssymb,euscript,amsthm,amsfonts}
\newtheorem{lem}{Lemma}
\newtheorem{thm}{Theorem}
\theoremstyle{definition}

\theoremstyle{remark}

\begin{document}

\renewcommand{\proofname}{Proof}
\makeatletter \headsep 10 mm \footskip 10 mm
\renewcommand{\@evenhead}%
{\vbox{\hbox to\textwidth{\strut \centerline{{\it Mathieu Dutour and Frank Vallentin}}} \hrule}}

\renewcommand{\@oddhead}%
{\vbox{\hbox to\textwidth{\strut \centerline{{\it
Some Six-Dimensional Rigid Forms
}}} \hrule}}

\begin{center}
{\Large\bf Some Six-Dimensional Rigid Forms}
\end{center}

\begin{center}
Mathieu DUTOUR (France)\footnote{Research financed by EC's IHRP Programme, within the Research Training Network ``Algebraic Combinatorics in Europe,'' grant HPRN-CT-2001-00272.} and Frank VALLENTIN (Germany)\footnote{Partially supported by the Edmund Landau Center for Research in Mathematical Analysis and Related Areas, sponsored by the Minerva Foundation (Germany).}
\end{center}

\vskip 15pt

\begin{quote}
One can always decompose Dirichlet-Voronoi polytopes of lattices
non-trivially into a Minkowski sum of Dirichlet-Voronoi polytopes of
rigid lattices. In this report we show how one can enumerate all rigid
positive semidefinite quadratic forms (and thereby rigid lattices) of a
given dimension~$d$. By this method we found all rigid positive
semidefinite quadratic forms for $d = 5$ confirming the list of $7$
rigid lattices by Baranovskii and Grishukhin. Furthermore, we found
out that for $d \leq 5$ the adjacency graph of primitive $L$-type
domains is an infinite tree on which $\mathsf{GL}_d(\mathbb{Z})$
acts. On the other hand, we demonstrate that in $d = 6$ we face a
combinatorial explosion.
\end{quote}

\section{Introduction}

Let $L$ be a lattice in Euclidean space~$(\mathbb{R}^d, \langle \cdot,
\cdot \rangle)$. With $L$ one associates the \textit{Dirichlet-Voronoi
polytope}
$$
\mathrm{DV}(L) = \{\mathbf{x} \in \mathbb{R}^d : 
\mbox{$\langle \mathbf{x}, \mathbf{x} \rangle \leq \langle \mathbf{v} - \mathbf{x}, \mathbf{v} - \mathbf{x}\rangle$ for all $\mathbf{v} \in L$}\}.
$$ 
One can always decompose Dirichlet-Voronoi polytopes of lattices
non-trivially into a Minkowski sum of Dirichlet-Voronoi polytopes of
so-called rigid lattices. In this paper we show how one can find all
rigid lattices in a given dimension.

It will be more convenient to use the language of quadratic
forms. With $L$ we associate a positive definite quadratic form: By
choosing a lattice basis $\mathbf{b}_1, \ldots, \mathbf{b}_d$ we get a
positive definite matrix $Q = (\langle \mathbf{b}_i,
\mathbf{b}_j\rangle)_{i,j}$ giving the positive definite quadratic
form $\mathbf{x} \mapsto \mathbf{x}^t Q \mathbf{x}$.

By applying basic facts of Voronoi's theory of $L$-type domains we get
an algorithm for finding all rigid forms of a given dimension. All
quadratic forms whose Dirichlet-Voronoi polytopes have the same
combinatorial-metric structure belong to a pointed polyhedral cone, a
so-called $L$-type domain. $L$-types domains of maximal dimension are
called \textit{primitive}, $L$-type domains of minimal dimension $1$
are called \textit{rigid}, and elements of rigid $L$-type domains are
called \textit{rigid forms}.  From Voronoi's algorithm for finding all
primitive $L$-type domains we get the facets of every primitive
$L$-type domain. By converting the facet description we find all
extreme rays. If an $L$-type domain is invariant under a non-trivial
symmetry group we can speed up the task of converting considerably
using the adjacency decomposition method.

We applied this algorithm for dimensions $\leq 5$. Thereby, we
confirmed the results by Baranovskii and Grishukhin (2): There is
exactly one $1$-dimensional rigid form (the Dirichlet-Voronoi polytope
is a line segment), there are no $2$- and $3$-dimensional rigid forms,
there is exactly one $4$-dimensional rigid form (the Dirichlet-Voronoi
polytope is the $24$-cell), and there are seven $5$-dimensional rigid
forms. Using this classification we verified that for $d \leq 5$ the
graph of primitive $L$-type domains is an infinite tree on which the
group $\mathsf{GL}_d(\mathbb{Z})$ acts.  It is computationally quite
simple to perform these classifications because face lattices of these
primitive $L$-type domains are very tame.

\textit{This is no longer the case in dimension~$6$.} We explored two
notable primitive $L$-type domains of $6$-dimensional quadratic
forms. The first cone has $130$ facets and we do not know a primitive
$L$-type domain in dimension~$6$ having more facets. The second cone
has $100$ facets and it contains a positive definite quadratic form
associated to the lattice whose covering density is conjectured to be
optimal in dimension~$6$. The automorphism groups of both cones are
fairly big. By using the adjacency decomposition method we succeeded
to compute the extreme rays of both cones: The fist one has
$7,145,429$ extreme rays, and the second one has $2,257,616$ extreme
rays. But many of these extreme rays correspond to equivalent (under
the group $\mathsf{GL}_d(\mathbb{Z})$) rigid forms. In total we found $25,263$
non-equivalent rigid $6$-dimensional positive definite quadratic
forms.

This paper is organized as follows: In Section~\ref{sec:ltype} we
recall some definitions and facts of Voronoi's theory of $L$-type
domains. In Section~\ref{sec:geometry} we show in which sense rigid
forms are building blocks of Dirichlet-Voronoi polytopes.  In
Section~\ref{sec:adjacency} we explain the adjacency decomposition
method.  In Section~\ref{sec:dim5} and Section~\ref{sec:dim6} we
report on computational results for rigid forms up to dimension~$6$.

\section{Notation: $L$-Type Domains}
\label{sec:ltype}

In this section we define $L$-type domains for positive semidefinite
quadratic forms. This enables us to define rigid forms.

By $\mathcal{S}^d$ we denote the space of all quadratic forms in $d$
variables, and by $\mathcal{S}^d_{\geq 0}$ we denote the set of all
positive semidefinite quadratic forms which is a closed pointed
cone. In the following we will identify $\mathcal{S}^d$ with the space
of all symmetric $(d \times d)$-matrices. We say that two quadratic
forms $Q, Q'$ are \textit{arithmetically equivalent} if there exists
an integral unimodular matrix $A \in \mathsf{GL}_d(\mathbb{Z})$ so
that $Q' = A^t Q A$.

Let $Q \in \mathcal{S}^d_{\geq 0}$ be a positive semidefinite
quadratic form arithmetically equivalent to $\left(\begin{smallmatrix}
Q' & 0\\ 0 & 0\end{smallmatrix}\right)$ where $Q'$ is positive
definite.  We define the \textit{Dirichlet-Voronoi polytope} of $Q$ by
$$
\mathrm{DV}(Q) = \{\mathbf{x}^t Q \in ({\mathbb R}^d)^* :
\mbox{$\mathbf{x}^t Q \mathbf{x} \leq (\mathbf{x}-\mathbf{v})^t Q (\mathbf{x}-\mathbf{v})$ for all $\mathbf{v} \in \mathbb{Z}^d$}\}.
$$

This way of defining Dirichlet-Voronoi polytopes of positive semidefinite quadratic forms is due to Namikawa (7). The definition has the important feature that we can define $L$-type domains of positive semidefinite quadratic forms by using the concept of strongly isomorphic polytopes which we recall now.

Let $V$ be a finite dimensional real vector space. Let $P \subseteq V$
be a convex polytope, and let $f \in V^*$ be a linear functional. We
define the \textit{support functional} of $P$ by $\eta(P, f) =
\max\{f(x) : x \in P\}$ and by $P_f = \{x \in P : f(x) = \eta(P, f)\}$
we denote the face of $P$ in direction $f$. We say that two convex
polytopes $P, Q \in V$ are \textit{strongly isomorphic} if for every
$f, g \in V^*$ with $P_f \subseteq P_g$ we have $Q_f \subseteq Q_g$.

We say that two positive semidefinite quadratic forms $Q, Q'$ belong
to the same \textit{$L$-type} if their Dirichlet-Voronoi polytopes are
strongly isomorphic. The set of all positive semidefinite quadratic
forms belonging to the same $L$-type is called an \textit{$L$-type
domain}.  Two $L$-type domains $\Delta$, $\Delta'$ are
\textit{arithmetically equivalent} if there exists $A \in
\mathsf{GL}_d(\mathbb{Z})$ so that $\Delta' = A^t \Delta A$.

In (11) Voronoi showed that $L$-type domains are open pointed
polyhedral cones, that the $L$-type domains give a face-to-face
partition of $\mathcal{S}^d_{\geq 0}$, and that there are only
finitely many non-equivalent $L$-type domains. $L$-type domains which
are of maximal dimension $\frac{d(d+1)}{2}$ are called
\textit{primitive}. $L$-type domains which are of minimal dimension
$1$ are called \textit{rigid}. Positive semidefinite quadratic forms
lying in a rigid $L$-type domain are called \textit{rigid}, too.

\section{Rigid Forms and Dirichlet-Voronoi Polytopes}
\label{sec:geometry}

In this section we show in what sense Dirichlet-Voronoi polytopes of
rigid positive semidefinite quadratic forms are building blocks of
Dirichlet-Voronoi polytopes of general positive semidefinite quadratic
forms.

\begin{lem}
\label{lem:minkowskivoronoi}
Every Dirichlet-Voronoi polytope of a positive semidefinite quadratic
form is Minkowski sum of Dirichlet-Voronoi polytopes of rigid
forms. More precisely: Let~$\overline{\Delta}$ be the topological
closure of an $L$-type domain. For positive semidefinite quadratic
forms $Q_1, \ldots, Q_n \in \overline{\Delta}$ and non-negative
numbers $\alpha_1, \ldots, \alpha_n$ we have
$$
\mathrm{DV}(\sum_{i=1}^n \alpha_i Q_i) = \sum_{i=1}^n \alpha_i \mathrm{DV}(Q_i).
$$
\end{lem}

The authors do not know exactly the origin of this lemma. Loesch gave
it in a dual formulation in (6). Later, Ryshkov gave in (9) a similar
but less precise statement.

\section{Computational techniques}
\label{sec:adjacency}

Fukuda's program \texttt{cdd} (5) computes the list of extreme rays of
a polyhedral cone given its list of facets. In our case the number of
extreme rays can be very large so that we cannot use \texttt{cdd}
naively. We apply another technique called adjacency decomposition
method to use the symmetry of the polyhedral cones we are considering.

Let $C \subseteq \mathbb{R}^d$ be $d$-dimensional polyhedral cone
determined by a set of facets $\{F_1, \ldots, F_n\}$. We assume that
$C$ is pointed at the origin. By $f_i \in (\mathbb{R}^d)^*$ we denote
a linear functional defining $F_i$, i.e.\ $F_i \subseteq \{\mathbf{x}
\in \mathbb{R}^d : f_i(\mathbf{x}) = 0\}$. Let $E$ be an initial
extreme ray of $C$ which we find e.g.\ by solving a generic linear
program on $C$. We compute the extreme rays adjacent to $E$: first we
project $C$ along $E$ by a linear map $\pi$.  Then, we find the
extreme rays of this projected cone $\pi(C)$ (using \texttt{cdd}, or
applying this procedure recursively). Every extreme ray $E_{\pi}$ of
$\pi(C)$ corresponds to a two-dimensional face $F$ of $C$ in which $E$
lies. Therefore, there is exactly one more extreme ray $E'$ of $C$ in
$F$. Every $\mathbf{e}' \in E'$ can we written as $\mathbf{e}' =
\alpha \mathbf{e} + \beta \mathbf{e}_{\pi}$, with $\mathbf{e} \in E$,
$\mathbf{e}_{\pi} \in E_{\pi}$ and some $\alpha, \beta$, which we have
to compute. This can be done by solving a ``two-dimensional linear
program": $f_i(\alpha \mathbf{e} + \beta \mathbf{e}_{\pi}) = 0$ for
all facets of $C$ incident to $F$, $f_j(\alpha \mathbf{e} + \beta
\mathbf{e}_{\pi}) \geq 0$ for all other facets of $C$. The key
computational step in the procedure above is the computation of the
extreme rays of $\pi(C)$. The complexity of this computation is
related to the incidence number of $E$, i.e.\ the number of facets
containing $E$.

The adjacency decomposition method applies to polyhedral cones having
a non-trivial symmetry group:
\begin{enumerate}
\item Take an initial list of orbits of extreme ray of $C$.
\item Take a representative $E$ of an orbit and finds the extreme rays $(E_i)_{1\leq i\leq m}$ adjacent to it.
\item If some $E_i$ represents a new orbit, then we add it to the list of orbits.
\item Finish when all orbits have been treated.
\end{enumerate}
This procedure has two main computational bottlenecks: it can be
difficult to identify new orbits, and the incidence of extreme rays
can be too high.

Since the symmetry groups of the cones we considered was not too big,
the first bottleneck was not a problem: for every new extreme ray we
generated the whole orbit.

For dealing with the second bottleneck we used Balinski's theorem:

\begin{thm} ((1), see e.g.\ (12))\\
Let $C$ be a $d$-dimensional pointed polyhedral cone. Let $G$ be the
undirected graph whose vertices are the extreme rays of $C$ and whose
edges are the $2$-dimensional faces of $C$. Two vertices $E_1, E_2$
are connected by an edge $F$ if $E_1, E_2 \in F$. Then, the graph $G$
is $(d - 1)$-connected, i.e.\ removal of any $d - 2$ vertices
leaves it connected.
\end{thm}

Due to Balinski's theorem, we can replace the criterion ``Finish when
all orbits have been treated'' by ``Finish when the number of extreme
rays in untreated orbits is lower than $d-2$''.

We applied the adjacency decomposition technique, starting with orbits
of lowest incidence. After some time, we found the complete list of
orbits. But still we had to treat the orbits with highest
incidence. These very degenerate orbits have usually a particular
significance, in our case they correspond to quadratic forms lying in
the boundary of $\mathcal{S}^d_{\geq 0}$. In the $L$-type domain
considered, we find out a peculiarity: the number of elements of the
orbits with highest incidence is very low, actually lower than $d - 2$.
Therefore, those orbits cannot disconnect the skeleton graph and so we
can stop earlier avoiding the computation of adjacencies of those
orbits.

\section{Dimensions $1, \ldots, 5$.}
\label{sec:dim5}

In dimension $1$ all positive definite quadratic forms are rigid;
their Dirichlet-Voronoi polytopes are line segments.  In dimension $2$
and $3$ there are no rigid positive definite quadratic forms because
in these dimensions every Dirichlet-Voronoi polytope is a zonotope
whence it is a Minkowski sum of line segments. In all these dimensions
there is only one non-equivalent primitive $L$-type domain.  In
dimension $4$ there are three non-equivalent primitive $L$-type
domains, and there is exactly one rigid positive definite quadratic
form which is associated to the root lattice $\mathsf{D}_4$. Its
Dirichlet-Voronoi polytope is the $24$-cell.  In dimension $5$ there
are $7$ rigid positive definite quadratic forms. First, they were
enumerated by Baranovskii and Grishukhin (2) by using Engel's list of
zone-contracted lattices (4). Our computations confirmed their result.

Now we will argue that these computations do not require much
computational effort. All $L$-type domains up to dimension~$4$ are
simplicial polyhedral cones. As noticed first by Barnes and Trennery
(3) (see also the discussion in (8) \S 13) this does no longer hold in
dimension~$5$ and above.  The following table shows how the numbers of
facets of primitive $L$-type domains are distributed among the $222$
non-equivalent primitive $L$-types domains in dimension~$5$. With $n$
we denote the number of facets and with $L_1(n)$ we denote the number
of non-equivalent primitive $L$-type domains in dimension~$5$ having
exactly $n$~facets.

\begin{center}
\begin{tabular}{r||c|c|c|c|c|c|c|c|c|c|c|c|c}
$n$ & $15$ & $16$ & $17$ & $18$ & $19$ & $20$ & $21$ & $22$ & $23$ & $24$ & $25$& $26$ & $27$\\
\hline
$L_1(n)$ & $62$ & $61$ & $46$ & $17$ & $10$ & $15$ &  $6$ &  $0$ &  $1$ &  $3$ & $0$ &  $0$ &  $1$\\
\end{tabular}
\end{center}

The next table shows the distribution of the numbers of extreme rays
among the $222$ non-equivalent primitive $L$-type domains in
dimension~$5$. With~$n$ we denote the number of extreme rays and with
$L_2(n)$ we denote the number of non-equivalent primitive $L$-type
domains having exactly $n$ extreme rays.

\begin{center}
\begin{tabular}{r||c|c|c|c|c|c|c|c|c|c|c|c}
$n$    & $15$ & $16$ & $17$ & $18$ & $19$ & $20$ & $21$ & $22$ & $23$ & $24$ & $25$ & $26$\\
\hline
$L_2(n)$ & $62$ & $84$ & $13$ &  $5$ & $33$ & $13$ &  $6$ &  $0$ &  $0$ &  $0$ & $0$ & $6$\\
\end{tabular}
\end{center}

The following table shows how the ranks of the extreme ray are
distributed among the $222$ non-equivalent $L$-type domains in
dimension $5$.  By $n$ we denote the number of extreme rays of and by
$R_k(n)$ we denote the number of non-equivalent $L$-type domains
having exactly $n$ extreme rays containing positive semidefinite
quadratic forms of rank~$k$, $k \in \{1,4,5\}$.

\begin{center}
{\footnotesize
\begin{tabular}{r||c|c|c|c|c|c|c|c|c|c|c|c|c|c|c|c|c}
$n$ & $0$ & $1$ & $2$ & $3$ & $4$ & $5$ & $6$ & $7$ & $8$ & $9$ & $10$& $11$ & $12$ & $13$ & $14$ & $15$ & $16$\\
\hline
\hline
$R_1(n)$ & $0$ & $0$ & $0$ & $0$ & $0$ & $0$ & $0$ & $0$ & $0$ & $0$ & $135$ & $58$ & $24$ & $3$ & $1$ & $1$ & $0$\\
\hline
$R_4(n)$ & $55$ & $49$ & $92$ & $19$ & $0$ & $7$ & $0$ & $0$ & $0$ & $0$ & $0$ & $0$ & $0$ & $0$ & $0$ & $0$ & $0$\\
\hline
$R_5(n)$ & $2$ & $12$ & $29$ & $38$ & $56$ & $14$ & $13$ & $17$ & $17$ & $8$ & $4$ & $0$ & $6$ & $0$ & $0$ & $0$ & $6$\\
\hline
\end{tabular}
}
\end{center}

By our computations we found out that for $d \leq 5$ the graph of
primitive $L$-type domains (vertices = primitive $L$-type domains,
edges = facets between primitive $L$-type domains) is an infinite tree
on which the group $\mathsf{GL}_d(\mathbb{Z})$ acts. The graph of
primitive L-type domains is a tree if and only if it has no cycle and
therefore if and only if every ridge contains at least one degenerate
form. Is this always the case? We think that the answer is ``No'' even
for $d = 6$. We also think that there is a primitive L-type domain
whose extreme rays are all non-degenerate.

\section{Dimension $6$}
\label{sec:dim6}

We consider here two different primitive $L$-type of $6$-dimensional
positive semidefinite quadratic forms, which were considered by the
second author in (10).

\subsection{The Cone $C_1$}

The cone $C_1$ is a primitive $L$-type domain of $6$-dimensional
positive semidefinite quadratic forms. It has $130$ facets and we
conjecture that there is no $L$-type domain of $6$-dimensional
positive semidefinite quadratic forms having more facets.  The
automorphism group of $C_1$ has order $1920$.  Using the adjacency
decomposition method we computed that $C_1$ has $7,145,429$ extreme
rays in $4,440$ orbits. It is worthwhile to note that among the
$4,440$ orbits there are two orbits which are equivalent under the
group $\mathsf{GL}_6(\mathbb{Z})$, and that there is one orbit which does
contain forms of rank~$5$, so that we found altogether
$4,438$ non-equivalent (under the group $\mathsf{GL}_6(\mathbb{Z})$) rigid
positive definite quadratic forms.

\subsection{The Cone $C_2$}

The cone $C_2$ contains the $6$-dimensional positive definite
quadratic form associated to the best known $6$-dimensional lattice
covering (10). It has $100$ facets and its automorphism group has
order $120$. Using the adjacency decomposition method we computed that
$C_2$ has $2,257,616$ extreme rays in $20,871$ orbits which correspond
to $20,861$ non-equivalent (under the group $\mathsf{GL}_6(\mathbb{Z})$) rigid
positive definite quadratic forms. Three orbits correspond to positive
semidefinite quadratic forms being not positive definite.

\subsection{Connection between $C_1$ and $C_2$}

In this section we show how $C_1$ and $C_2$ are related.
Both cones $C_1$ and $C_2$ contain the rigid form
$$
Q_{\mathsf{E}_6^*} =
\begin{pmatrix}
 4 & 1 & 2 & 2 &-1 & 1\\
 1 & 4 & 2 & 2 & 2 & 1\\
 2 & 2 & 4 & 1 & 1 & 2\\
 2 & 2 & 1 & 4 & 1 & 2\\
-1 & 2 & 1 & 1 & 4 & 2\\
 1 & 1 & 2 & 2 & 2 & 4
\end{pmatrix}
$$
associated to the lattice $\mathsf{E}_6^*$. The automorphism groups of
$G_i$ of $C_i$, $i = 1,2$ are subgroups of the automorphism group $G
= \{T \in \mathsf{GL}_6(\mathbb{Z}) : T^t Q_{\mathsf{E}_6^*} T =
Q_{\mathsf{E}_6^*}\}$ of $Q_{\mathsf{E}_6^*}$. 

The subspace $\mathbf{I}_1$ of all quadratic forms invariant under the
group $G_1$ is spanned by $Q_{\mathsf{E}^*_6}$ and $R_1$ (see below)
which is an extreme ray of $C_1$. If we intersect the cone $C_2$ with
$\mathbf{I}_1$ we get a two-dimensional cone with extreme rays
$Q_{\mathsf{E}^*_6}$ and $R_2$. The rigid forms $R_1$ and $R_2$ are
$$
R_1 = 
\begin{pmatrix}
12 & 3 & 6 & 6 & -3 & 3\\
 3 & 7 & 4 & 4 &  3 & 2\\ 
 6 & 4 & 8 & 3 &  1 & 4\\
 6 & 4 & 3 & 8 &  1 & 4\\ 
-3 & 3 & 1 & 1 &  7 & 3\\
 3 & 2 & 4 & 4 &  3 & 7
\end{pmatrix}
\quad
R_2 =
\begin{pmatrix}
0 & 0 & 0 & 0 & 0 & 0\\
0 & 5 & 2 & 2 & 3 & 1\\
0 & 2 & 4 & 0 & 2 & 2\\
0 & 2 & 0 & 4 & 2 & 2\\
0 & 3 & 2 & 2 & 5 & 3\\
0 & 1 & 2 & 2 & 3 & 5
\end{pmatrix}
$$

\subsection{Further Remarks}

It is remarkable that there are only very few instances of
arithmetically equivalent extreme rays, which are not equivalent under
the symmetry group of the cone. The number of extreme rays of both
cones is extremely large. Nevertheless it is interesting to note that
the number of non-equivalent extreme rays corresponding to forms which
are not positive definite is low. In total, we obtained $25,263$ new
rigid positive definite quadratic forms in dimension~$6$.

\vskip 10 pt

\centerline{\Large REFERENCES} \vskip 7 pt

\begin{enumerate}
\small
\item M.L.\ Balinski (1961), {On the graph structure of convex polyhedra in $n$-space}, Pacific J. Math., {\bf 11}, 431--434
\item E.P.\ Baranovskii, V.P.\ Grishukhin (2001), {\it Non-rigidity degree of a lattice and rigid lattices}, European J. Combin., {\bf 22}, 921--935.
\item E.S.\ Barnes, D.W.\ Trenerry (1972), {\it A class of extreme lattice-coverings of $n$-space by spheres}, J.\ Austral.\ Math.\ Soc., {\bf 14}, 247--256.
\item P.\ Engel (1998), {\it Investigations of parallelohedra in $\mathbb{R}^d$}, in {\sl Voronoi's impact on modern science}. Proc.\ Math.\ Nat.\ Acad.\ Sci.\ Ukraine, \textbf{21}, 22--60.
\item K. Fukuda (2001), {\it cddlib reference manual, cddlib Version 0.92}, ETH Z\"urich, {\tt http://www.ifor.math.ethz.ch/\~{}fukuda/cdd home/cdd.html}
\item H.-F.\ Loesch (1990), {\it Zur Reduktionstheorie von Delone-Voronoi f\"ur matroid\-ische quadratische Formen}, PhD Thesis, Ruhr-Universit\"at Bochum.
\item Y.\ Namikawa (1976), {\it A new compactification of the Siegel space and degenerations of abelian varieties, I, II}, Math.\ Ann., {\bf 221}, 97--141, and 201--241.
\item S.S.\ Ryshkov, E.P.\ Baranovskii (1976), {\it C-types of $n$-dimensional lattices and 5-dimensional primitive parallelohedra (with application to the theory of coverings)}, Proc.\ Steklov Inst.\ Math., {\bf 137}, 1--140.
\item S.S.\ Ryshkov (1998), {\it On the structure of a primitive parallelohedron and Voronoi's last problem}, Russian Math.\ Surveys, {\bf 53}, 403--405.
\item F.\ Vallentin (2003), {\it Sphere Coverings, Lattices, and Tilings (in Low Dimensions)}, PhD Thesis, Munich University of Technology, Online Publication:\\
{\footnotesize
{\tt http://tumb1.biblio.tu-muenchen.de/publ/diss/ma/2003/vallentin.html}
}
\item G.F.\ Voronoi (1909), {\it Nouvelles applications des parame\'etres continus \`a l\`a th\'eorie des formes quadratiques, Deuxi\`eme M\'emoire,
Recherches sur les parall\'elloedres primitifs}, J.\ Reine Angew.\ Math.\ {\bf 134}, 198--287 and {\bf 136}, 67--181.
\item G.M.\ Ziegler (1995), {\it Lectures on polytopes}, New-York: Springer-Verlag.

\end{enumerate}

\vskip 20 pt

\noindent {\it The Hebrew University of Jerusalem, Israel
}

\vskip 2 pt

\noindent {\it e-mail: {\tt mathieu.dutour@ens.fr} and {\tt vallenti@ma.tum.de}}

\end{document}